\documentclass{amsart}

\usepackage{fancyhdr}
\usepackage{amssymb}
\usepackage{amsmath}
\usepackage{latexsym}
\input xy 
\xyoption{all}

\newtheorem{teo}{Theorem}[section]
\newtheorem{prop}[teo]{Proposition}
\newtheorem{lem}[teo]{Lemma}
\newtheorem{cor}[teo]{Corollary}

\newtheorem{defi}[teo]{Definition}
\newtheorem{fac}[teo]{Fact}
\newtheorem{rem}[teo]{Remark}

\def\dnfo{\,\raise.2em\hbox{$\,\mathrel|\kern-.9em\lower.35em\hbox{$\smile$}$}}
\def\dfo{\;\raise.2em\hbox{$\mathrel|\kern-.9em\lower.35em\hbox{$\smile$}\kern-.7em\hbox{\char'57}$}\;}
\newcommand{\s}{SU}

\newcommand{\si}{\sigma}
\newcommand{\na}{\mathbb{N}}

\newcommand{\gr}{\mathbb G}
\newcommand{\sse}{\subseteq}
\begin{document}
\title{Definable Groups in DCFA}

\author{Ronald F. Bustamante Medina}
\address{Centro de Investigaciones en Matem\'atica Pura y Aplicada\\
Escuela de Matem\'aticas\\
Universidad de Costa Rica\\
Sede Rodrigo Facio,
2060, San Jos\'e, Costa Rica}
\email{ronald.bustamante@ucr.ac.cr}
\date{April 12, 2018}

\begin{abstract}
E. Hrushovski proved that the theory of difference-differential fields of characteristic zero has a model-companion. We denote it $\it DCFA$.
In this paper we study definable groups in a model of $\it DCFA$. First we prove that such a group is embeds on an algebraic group. Then we study 
1-basedeness, stability and stable embeddability of abelian definable groups.
\end{abstract}

\maketitle

\section{Preliminaries}
The class of differentially closed fields of characteristic zero with a generic automorphism is elementary, we denote it {\it DCFA}.   

Our aim in this paper is to study definable groups in models of {\it DCFA}: in section~\ref{sec:embed}  we prove that a definable group in a model of {\it DCFA} embeds in an algebraic group. In section~\ref{sec:stable} prove that we can reduce questions about 1-basedness and stable, stable embeddability in {\it DCFA} to questions about 1-basedness and stable, stable embeddability in either {\it DCF} or {\it ACFA}. We use this in section~\ref{sec:abelian} to study the model theory of definable abelian groups.

We give now a brief summary of what we know about {\it DCFA}. Since we will work in difference, differential and difference-differential fields we will denote
the respective languages by $\mathcal{L}_{\si}$, $\mathcal{L}_{D}$ and $\mathcal{L}_{\si, D}$

In \cite{rbdcfa1} we give an axiomatisation of {\it DCFA} and prove its main properties: given a model of {\it DCFA} it is of course a differentially closed field (model of {\it DCF}) and an algebraically closed field with a generic automorphism (model of {\it ACFA}). Independence is defined by linear disjointness.  This theory is not complete, but its completions are
easily described,  those completions eliminate imaginaries (moreover, they satisfy the Independence
Theorem over algebraically closes sets) and thus are supersimple and types are ranked by
the $\s$-rank. Forking is determined by quantifier-free formulas, thus {\it DCFA} is
 quantifier-free $\omega$-stable.
A basis theorem for  (perfect) difference-differential ideals imply that in a model of {\it DCFA} the difference-differential Zariski topology (defined in analogy with Zariski topology
in algebraically closed fields) is Noetherian.

Let $(K, \si, D)$ be  a model of {\it DCFA}, there are two important definable subfields of $K$, the field of constants ${\mathcal C}=\{x \in K: Dx=0\}$ and the fixed field $Fix(\si)=\{x \in K: \si(x)=x\}$. 

Given $a \in K$ and $A \sse K$, we define the $(\si,D)$-transcendence degree of $a$ over $A$
as the transcendence degree of the difference-differential field generated by $A$ and $a$ over $A$. In the cases of {\it DCF} and {\it ACFA} the finiteness of such a degree is equivalent to the finiteness of the rank of $a$ over $A$. However this does not hold for {\it DCFA}:
in \cite{rbrank1} we give an example of a set whose generic type has infinite
$(\si,D)$-transcendence degree but $\s$-rank 1). This represents a difficulty in the treatment of definable groups,
so we shall try different ways to describe  definable groups departing from properties of groups definable in
differential and difference fields.

In \cite{rbjets} and \cite{rbarcs01} we proved that Zilber's dichotomy holds for 
{\it DCFA}: a type of $\s$-rank 1 either has  a simple geometry (it is 1-based) or has a strong
interaction with (is non-orthogonal to)  $Fix(\si) \cap {\mathcal C}$.

We now introduce some definitions and useful facts about  definable groups in supersimple theories.
Let $T$  be a supersimple theory, $M$ a saturated model of $T$, let $G$ be a type-definable (definable by an infinite number of formulas) group
% and $G$ an
%type-definable (definable by an infinite number of formulas) group over some 
and let $A \subset M$ be a set of parameters.

\begin{defi}
Let $G$ be a 
Let $p \in S(A) $. We say that $p$ is a left generic type of $G$ over $A$ if it is realized in $G$ and 
for every $a \in G$ and $b$ realizing $p$ such that $a \dnfo_A b  $, we have $b \cdot a \dnfo_A a$.
\end{defi}

%Some of the properties of generic types in $\omega$-stable groups hold in simple theories.

The following result is proved in \cite{simpg} :

\begin{fac}%Let $G$ be an $A$-definable group.
\begin{enumerate}

\item Let $a,b \in G $. If $tp(a/Ab)  $ is left generic of $G$, then so is $tp(b \cdot a/Ab ) $.

\item Let $p \in S(A)  $ be realized in $G$, $B=acl(B) \supset A$, and $q \in S(B) $
 a non-forking extension of $p$. Then $p$ is a generic of $G$ if and only if $q$ is a generic of $G$.

\item Let $tp(a/A)  $ be generic of $G$; then so is $tp(a^{-1}/A)$.

\item There exists a generic type of $G$.

\item A type is left generic if and only if it is right generic.

\end{enumerate}
\end{fac}
The following fact is
proved in \cite{wag}, chapter 5.

\begin{fac}\label{GIII5}
Let  $H$ a type-definable subgroup of $G$,
\begin{enumerate}
\item Let $p \in S(A)  $, then $ p$ is a generic of $G$ over $A$ if and only if $\s(G)=\s(p) $.
\item $\s(G)=\s(H)$ if and only if $[H:G]< \infty  $.
\item $\s(H)+\s(G/H) \leq SU(G) \leq \s(H) \oplus \s(G/H) . $ 
\end{enumerate}
\end{fac}

\section{Every Definable Group Embeds in an Algebraic Group}\label{sec:embed}

We introduce $*$-definable groups in stable theories.
Suppose that $T$ is a complete theory and $M$ a saturated model of $T$.
A $*$-tuple is a tuple $(a_i)_{i \in I}$, where $I$ is an index set of cardinality less than the cardinality of $M$,
and $a_i \in M^{eq}$ for all $i \in I$.
Let $A \subset M$. A $*$-definable set is a collection of $*$-tuples, indexed by the same set of parameters $I$,
which is the set of realizations of  a partial type $p(x_i)_{i \in I}$ over $A$.
A $*$-definable group is a group with $*$-definable domain and multiplication.

The following propositions are proved in \cite{kopi}. Recall that the canonical base of a
strong  type $p$, $Cb(p)$ is the set that is fixed pointwise by the automorphisms that fix $p$.
\begin{prop}\label{GIII1}
Let $T$ be a stable theory; $M$ a saturated model of $T$.
Let $a,b,c,x,y,z$ be $*$-tuples of $M$ of length strictly less than the cardinal of $M$, such that:

\begin{enumerate}

\item $acl(M,a,b)=acl(M,a,c) = acl(M,b,c)       $

\item $acl(M,a,x)  =  acl(M,a,y)$ and $Cb(stp(x,y/M,a)) $ is interalgebraic with $a$ over $M $.

\item As in 2. with $b,z,y$ in place of $a,x,y $

\item As in 2. with $c,z,x$ in place of $a,x,y $

\item Other than $\{a,b,c  \},  \{a,x,y  \},  \{b,z,y  \},\{c,z,x  \} $, any 3-element subset of $\{a,b,c,x,y,z    \}  $ is independent over $\mathcal M$.
\end {enumerate}

Then there is a $*$-definable group $H$ defined over $M$ and $a',b',c'  \in H $ generic independent
 over $M$ such that $a$  is interalgebraic with $a'$ over $M$,   $b$  is interalgebraic with $b'$ over $M$
 and  $c$  is interalgebraic with $c'$ over $M$.

\end{prop}

\begin{prop}\label{GIII2}
Let $T$ be a simple theory; $M$ a saturated model of $T$. Let $G,H$ be type-definable groups,
 defined over $K \prec {M}$, and let $a,b,c \in G $  and $a',b',c' \in H  $ such that

\begin{enumerate}
\item $a,b$ are generic independent over $M$.
\item $a \cdot b = c  $ and $a' \cdot b' =c'  $.
\item $a$ is interalgebraic with $a'$ over $M$,  $b$ is interalgebraic with $b'$ over $M$ and
 $c$ is interalgebraic with $c'$ over $M$
\end{enumerate}
Then there is a type-definable over $M$ subgroup $G_1$ of bounded index in $G$,
and a type-definable over $M$ subgroup $H_1$ of $H$ and a  type-definable over $M$
isomorphism $f$ between $G_1/N_1  $ and $H_1 / N_2$ where $N_1  $ and $N_2$  are finite normal subgroups of $G_1 $ and $H_1$ respectively.

\end{prop}

\begin{rem}
If $T$ in \ref{GIII2} is supersimple and $G,H$ are definable,
then we can choose $G_1$ definable of finite index in $G$ and $f$ definable.
\end{rem}
The following result is proved in \cite{Hru}:
\begin{prop}\label{GIII3}
Let $G$ be a $*$-definable group in a stable structure.
Then there is a projective system of definable groups with inverse limit $G'$,
and a $*$-definable isomorphism between $G$ and $G'$.
\end{prop}

In \cite{difgp} the author proved that a ${\mathcal L}_D$-definable (definable in the language of 
differential fields) group in {\it DCF} is essentially a differential algebraic group and that a definable group in {\it DCF} virtually embeds in an algebraic group.

So, to prove that a definable group in {\it DCFA} embeds in an  algebraic group we will show that it embeds in a ${\mathcal L}_D$-definable group.

\begin{teo}\label{GIII6}
Let $({\mathcal U},\si,D)$ be a model of {\it DCFA}, $K \prec {\mathcal U} $ and $G$ a $K$-definable group.
Then there is an  ${\mathcal L}_D$-definable group $H$, a definable subgroup $G_1$  of $G$ of finite index,
and a definable isomorphism between $G_1/N_1$ and $H_1 / N_2$, where  $H_1$ is a definable subgroup of $H(\mathcal{U})$,
 $N_1$ is a finite normal subgroup of $G_1$, and  $N_2$ is a  finite normal subgroup of $H_1$.
\end{teo}

{\it Proof:}\\

Let $a,b,y   $ be generic independent elements of $G$ over $K$. Let $x = a \cdot y,  z = b^{-1} \cdot y,  c = a \cdot b  $, so $x = c \cdot z  $.
%Let $\bar{a}=(D^i\si^j (a): i  \in {\mathbb N},   j \in {\mathbb Z})  $,
Let $\bar{a}=(\si^i (a):    j \in {\mathbb Z})  $, and similarly for $\bar{b}, \bar{c}, \bar{x},\bar{y},\bar{z} $.
Then, as the model-theoretic algebraic closure of a set is the differential-field-theoretic algebraic closure of the set closed by $\si$, working in {\it DCF},
 $\bar{a},\bar{b},\bar{c},\bar{x},\bar{y},\bar{z} $
 satisfy the conditions of \ref{GIII1}.
Thus there is a $*$-${\mathcal L}_D$-definable group $H$ over $K$, and generic $K$-independent elements $a^*,b^*, c^* \in H    $ such that
 $\bar{a}  $ is interalgebraic with  $a^*  $ over $K$, $\bar{b}  $ is interalgebraic with  $b^*  $ over $K$,
 $\bar{c}  $ is interalgebraic with  $c^*  $ over $K$  and
 $ c^*=  a^* \cdot  b^*   $ (the interalgebraicity, independence and generics  in the sense of {\it DCF}).

Since  {\it DCF} is $\omega$-stable, by \ref{GIII3}, $H$ is the inverse limit of $H_i, i \in \omega $,
 where the $H_i$ are ${\mathcal L}_D$-definable groups.

Let $\pi_i : H \longrightarrow H_i  $ be the $i$-th canonical epimorphism.
Let $a_i = \pi_i (a^*) $, $b_i = \pi_i (b^*) $ and $c_i = \pi_i (c^*) $. Then $a^*   $ is
interalgebraic with $(a_i)_{i \in \omega } $ over $K$ ,
 $b^*   $ is interalgebraic with $(b_i)_{i \in \omega}  $ over $K$ and
$c^*   $ is interalgebraic with $(c_i)_{i \in \omega}  $ over $K$, all interalgebraicities in the
sense of {\it DCF}.

Since for $i < j$, $a_i \in K(a_j)$  , $b_i \in K(b_j)$  and  $c_i \in K(c_j)$,
there is $i \in \omega$ such that $a$ is interalgebraic with $a_i$ over $K$,   $b$ is interalgebraic
 with $b_i$ over $K$ and  t$c$ is interalgebraic with $c_i$ over $K$  in the sense of {\it DCFA}.
So we can apply  \ref{GIII2} to $a,b,c \in G$ and $a_i,b_i,c_i \in H_i$.\\
$\Box$

\begin{cor}
 Let $G$ be a definable group.
Then there is an  algebraic group $H$, a definable subgroup $G_1$  of $G$ of finite index,
and a definable isomorphism between $G_1/N_1$ and $H_1 / N_2$, where  $H_1$ is a definable subgroup of $H(\mathcal{U})$,
 $N_1$ is a finite normal subgroup of $G_1$, and  $N_2$ is a  finite normal subgroup of $H_1$
\end{cor}
\section{Stability, Stable Embeddability and 1-basedness}\label{sec:stable}

In this section we discuss how to apply results 
from \cite{salinas} to obtain similar results in models of {\it DCFA}.
We also give a criterion for 1-basedness in {\it DCFA}.

We begin with general definitions and facts on supersimple theories.

$T$ will denote a supersimple theory which eliminates imaginaries. Let $M$ be a
saturated model of $T$.

Let us recall that two types $p,q$ over $A \sse M$ are orthogonal, denoted $p \perp q$, if for every set $B \supseteq A$ and every realisations $a, b$ of $p$ and $q$ respectively, $a \dnfo_B b$.

%Def of stable, stably embedded (cf appendix 1 of [CH])
\begin{defi}  
\begin{enumerate}
\item Let $A \subset M$ and let $S$ be an $(\infty)$-definable set over $A$. We say that $S$ is 1-based if
for every tuple $a$ of $S$ and every $B \supseteq A$, $a$ and $B$ are independent over
$acl(Aa) \cap acl(B)$.
\item A type is 1-based if the set of its realizations is 1-based.
\end{enumerate}
\end{defi}

The following useful result is proved in \cite{wagbase}.

\begin{prop}\label{wgr}
\begin{enumerate}
\item The  union of $1$-based sets is $1$-based.
\item If $tp(a/A)$ and $tp(b/Aa)$ are $1$-based, so is $tp(a,b/A)$.
\end{enumerate}
\end{prop}

We introduce now stable, stably embedded types (also called fully stable types).

\begin{defi}
A (partial) type $p$ over a set $A$ is stable, stably embedded
 if whenever $a$ realizes $p$ and $B\supset A$, then $tp(a/B)$ is definable.
Equivalently, let $P$ denote the set of realizations of $p$. Then $p$ is 
stable, stably embedded if and only if for all set $S\cap P^n$ where $S$ 
is  definable, there is a set $S'$ definable with parameters from $P$ 
and such that $S'\cap P^n=S\cap P^n$.
\end{defi}

%[Note: if $p$ is complete, this is what Shelah calls a stable type]. \\

The following result is proved in the Appendix of \cite{salinas}: 

\begin{lem}\label{st0} 
If $tp(b/A)$ and $tp(a/Ab)$ are stable, stably embedded, so is $tp(a,b/A)$. 
\end{lem}

\begin{rem}\label{sserem}
In \cite{salinas}, a certain property of models of {\it ACFA} (called superficial stability) is 
isolated, and guarantees that certain types over algebraically closed 
sets are stationary, and therefore definable. It follows from model 
theoretic considerations that if for any algebraically closed set $B$ 
containing $A$, $tp(a/B)$ is stationary, then  $tp(a/A)$ will be 
stable, stably embedded.
\end{rem}
%Since we will work in difference, differential and difference-differential fields we will denote
%the respective languages by $\mathcal{L}_{\si}$, $\mathcal{L}_{D}$ and $\mathcal{L}_{\si, D}$
\begin{lem}\label{lemsse}
Let $(K,\si)$ be a model of {\it ACFA}, $A=acl_\si(A)\subset K$ and 
$a\in K$. Then $tp(a/A)$ is  stationary if and only if $tp(a/A)\perp 
(\si(x)=x)$, where $acl_\si$ denotes the model-theoretic algebraic closure in {\it ACFA}.

\end{lem}

{\it Proof}:\\

Indeed, write $SU(a/A)=\omega k+n$, and let $b\in
acl_\sigma(Aa)$ be such that $SU(b/A)=n$. Then $tp(b/A)\perp
(\sigma(x)=x)$ and, by Theorem 4.11 of [2], $tp(acl_\sigma(Ab)/A)$ is
stationary. If $c\in acl_\sigma(Aa)$ satisfies some 
non-trivial difference equation over $acl_\sigma(Ab)$ then 
$SU(c/Ab)<\omega$ 
and therefore $c\in acl_\sigma(Ab)$. Hence, by Theorem 5.3 of [3], 
$tp(a/acl_\sigma(Ab))$ is stationary, and therefore so is $tp(a/A)$.

For the converse, there are independent 
realizations $a_1,\cdots,a_n$ of $tp(a/A)$, and elements 
$b_1,\cdots,b_m\in Fix(\si)$ such that $(a_1,\cdots,a_n)$ and 
$(b_1,\cdots,b_m)$ are not independent over $A$. Looking at the field of 
definition of the algebraic 
locus of $(b_1,\cdots,b_m)$ over $acl_\si(A,a_1,\cdots,a_n)$, there is 
some $b\in Fix( 
\si)\cap acl_\si(A,a_1,\cdots,a_n)$, $b \not\in A$. Then $tp(b/A)$ is not 
stationary: if $c \in Fix(\si)$ is independent from $b$ over $A$, then 
$tp(b/A)$ has two distinct non-forking extensions to $Ac$, one in which 
$\sqrt{b+c}\in Fix(\si)$, the other in which $\sqrt{b+c}\not\in Fix(\si)$.
Hence $tp(a_1,\cdots,a_n/A)$ is not stationary, and neither is 
$tp(a/A)$.\\
$\Box$

%Let $(K,\si)$ be a model  of {\it ACFA} of characteristic $0$, $a$ a tuple in 
%$K$, and $A=acl(A)\subset K$. If $tp(a/A)$ is  orthogonal to the 
%formula $\si(x)=x$, then $tp(a/A)$ is stationary (see \cite{salinas}, and 
%\cite{salinas2}). It follows that if $tp(a/A)$ is hereditarily orthogonal to 
%$\si(x)=x$, then $tp(a/A)$ is stable stably embedded. In that 
%case, it will also be $1$-based.

It is important to note that stationarity alone does not imply 
stability: if $a$ is transformally transcendental over 
$A=acl_\si(A)$ (a is not the root of a non-zero $\si$-polynomial over $A$),
 then $tp_{ACFA}(a/A)$ is stationary, but it is not stable.
%We know that no completion of {\it DCFA} is stable. As in the case of 
%completions of {\it ACFA}, it turns out that certain definable sets, 
%endowed  with the structure induced by the ambient model, are stable 
 %stably embedded. 
These results can be used to give sufficient conditions on types in 
{\it DCFA} to be stationary, and stable, stably embedded.

\begin{prop}\label{st1}
Let $(K,\si,D)$ be a model of DCFA, let $A=acl(A)\subset K$, and 
$a$ a tuple in $K$.
\begin{enumerate}
\item Assume that $tp_{ACFA}(a,Da,D^2a,\cdots/A)\perp \si(x)=x$. Then 
$tp(a/A)$ is stationary.\label{st1_1}

\item Assume that for every $n$, every extension of $tp_{ACFA}(D^na/Aa\cdots D^{n-1}a)$ is 
 orthogonal to $(\si(x)=x)$. Then $tp(a/A)$ is stable, 
stably embedded. It is also 1-based. \label{st1_2}

\item If $tp(a/A)$ has an extension that is not orthogonal to $(\si(x)=x)$, then 
$tp(a/A)$ is not stable, stably embedded. 
\end{enumerate}
\end{prop}

{\it Proof}:\\

1. As $tp_{ACFA}(a,Da,D^2a,\cdots/A)\perp \si(x)=x$,  \ref{lemsse}  implies that $tp_{ACFA}(a,Da,D^2a,\cdots/A)$ 
is stationary. Since the $tp(a/A) $ is 
determined by $tp_{ACFA}(a,Da,D^2a,\cdots/A)$ , $tp(a/A)$ is stationary: 
 Let $b,c$ be two realizations of non-forking extensions of $tp(a/A)$ to a set $B=acl(B) \supset A$.
As $tp_{ACFA}(a,Da,D^2a,\cdots/A)$ is stationary we have that
 $tp_{ACFA}(b,Db,D^2b,\cdots/B) =tp_{ACFA}(c,Dc,D^2c,\cdots/B)$. 
 If $\varphi(x)$ is an
$\mathcal{L}_{\si,D}(B)$-formula satisfied by $b$, then there is a $\mathcal{L}_{\si}(B)$-formula
$\psi(x_0,\cdots,x_k)$ such that $\phi(b)=\psi(b,Db,\cdots,D^kb)$; so we have 
$\psi(b,Db,\cdots,D^kb)\in tp_{ACFA}(b,Db,D^2b,\cdots/B) =tp_{ACFA}(c,Dc,D^2c,\cdots/B)$. This implies that
$tp(b/B)=tp(c/B)$, and thus $tp(a/A)$ is stationary.\\
2. By \ref{lemsse} for all $n \in \na$ and for all $B \supset A$, $tp_{ACFA}(D^na/Ba\cdots D^{n-1}a)$ is stationary. Thus, by \ref{sserem},
for all $n$,  $tp_{ACFA}(D^na/Aa\cdots D^{n-1}a)$ is stable, stably embedded and 1-based.
By \ref{st0}  stability, stable embeddability is preserved by extensions, 
hence $tp_{ACFA}(a,Da,\cdots/A)$ is stable, stably embedded, and this implies that all extensions
to algebraically closed sets are stationary. As above, we deduce that all extensions of
$tp(a/A)$ to algebraically closed sets are stationary, hence $tp(a/A)$ is stable, stably embedded. 
By \ref{wgr} we have also that $tp_{ACFA}(a,Da,\cdots/A)$ is 1-based.
As  $tp(a/A) $ is 
determined by $tp_{ACFA}(a,Da,D^2a,\cdots/A)$ , $tp(a/A)$ is 1-based.
% By the definition of independence 
%in difference-differential fields and the fact that $acl(A,a)=acl_{DCF}(A,a,Da,\cdots)$ $tp(a/A)$ is 1-based:
%Let $A\subset B=acl(B) \subset C=acl(C)$ , and let   $b$ be tuple
 %of realizations of $tp(a/A)$.
% Let us suppose that $acl(Bb) \cap  C=B$.
 %By hypothesis $tp_{ACFA}(a,Da, \cdots/A)$ is
 %$1$-based,  therefore  $(b,Db, \cdots)$ is independent from $C$ over $B$ in {\it ACFA}.
% As $A=acl(A)$, we have that for all $n \in \na$, $tp_{DCF}(\si^n(b)/A)$ is
% $1$-based. By \ref{wgr} $tp_{DCF}(b,\si(b),\ldots, \si^n(b)/B)$ is also
% $1$-based.
 %Hence, $(b, Db,\ldots, D^nb)$ is {\it ACFA}-independent from $C$ over $B$,
 %for every $n \in  \na$.
 %Then for every finite subset $S$ of 
%$acl(Bb)$, $B(S)$ is linearly disjoint from $C$ over $B$ (that is 
%because 
%every such $S$ is such that $B(S)$ contained in $acl_{\si}(B,b,Db, 
%\cdots, 
%D^nb)$ for some $n$). Thus by definition of linear disjointness 
%$acl(Bb)$ is linearly disjoint from $C$ over $B$. So
%$b$ is {\it DCFA}-independent from $C$ over $B$.

3. If $tp(a/K)$ is not hereditarily orthogonal to $\si(x)=x$ then there is $B=acl(B) \supset A$
such that $tp(a/B) \not \perp \si(x)=x$. 
Then there are independent 
realizations $a_1,\cdots,a_n$ of $tp(a/B)$, and elements 
$b_1,\cdots,b_m\in Fix(\si)$ such that $(a_1,\cdots,a_n)$ and 
$(b_1,\cdots,b_m)$ are not independent over $B$.

If we look at the field of 
definition of the algebraic 
locus of $(b_1,\cdots,b_m)$ over $acl(A,a_1,\cdots,a_n)$, we can find  $b\in Fix 
(\si)\cap acl(A,a_1,\cdots,a_n)$, $b \not\in A$.  
Then $tp(b/A)$ is not stationary: 
Let $c \in Fix(\si)$ be independent from $b$ over $A$, then 
$tp(b/A)$ has two distinct non-forking extensions to $Ac$, one in which 
$\sqrt{b+c}\in Fix(\si)$, the other in which $\sqrt{b+c}\not\in Fix(\si)$.
Hence $tp(a_1,\cdots,a_n/A)$ is not stationary, and neither is 
$tp(a/A)$.\\
% If $tp(a/K)$ is not hereditarily orthogonal to $\si(x)=x$ then there is $B=acl(B) \supset A$
%such that $tp(a/B) \not \perp \si(x)=x$. Then there is $n \in \na$ such that 
%$tp_{ACFA}(a,Da,\cdots, D^na/B) \not\perp \si(x)=x$. By \ref{lemsse}, $tp_{ACFA}(a,Da, \cdots, D^na/A)$ is not
%stationary, then $tp(a/B)$ is not stationary, hence it is not stable stably embedded.\\
%We claim that $tp_{ACFA}(a,Da,D^a,\ldots /A) \not\perp (\si(x)=x)$. 
%(1) If $tp(a/A)$ is not stationary there is some $m$ such that $tp_{ACFA}(a,Da,\ldots,D^ma/A)$ is not stationary. 
%Let $p_1,\ldots,p_n$ be (${\mathcal L}_\si$-)types of $\s$-rank $1$ or $\omega$ (computed in {\it ACFA}) such that 
%$tp_{ACFA}(a,Da,\cdots,D^ma/A)$ is domination equivalent to $p_1\times \cdots \times p_n$.
% Then one of the $p_i$'s is not stationary, and by 5.3 of \cite{zh1}, it must have $\s$-rank $1$.
% Theorem 4.11 of \cite{salinas} finishes the proof.\\
%(2) is implied by (1) and \ref{wgr}.\\
$\Box$

\begin{rem}\label{st2}
Let $A,K$ and $a$ be as above.
\begin{enumerate}
\item If $\s(a/A)=1$, then the stationarity of $tp(a/A)$ implies its 
stability and stable embeddability.

\item There are examples of types of $\s$-rank $1$ which satisfy \ref{st1}(\ref{st1_1})
above but do not satisfy \ref{st1}(\ref{st1_2}). Thus condition \ref{st1}(\ref{st1_2}) is not 
implied by stationarity.

%\item $tp(ab/A)$ is stable stably embedded (resp. 1-based) if and 
%only if $tp(a/A)$ and $tp(b/Aa)$ are stably embedded (resp. 
%1-based). (See the appendix of \cite{salinas} and \cite{wag}).
\end{enumerate}
\end{rem}

\begin{cor}\label{st3}
Let $A=acl(A)$, and $a$ a tuple in ${\mathcal C}$.
 Then $tp(a/A)$ is stable, stably embedded if and only if 
$tp_{ACFA}(a/A)$ is stable, stably embedded. In this case, it will also be $1$-based. 
\end{cor}

\begin{prop}\label{st4}
Let $A=acl(A)\subset K$, and $a$ a tuple in $K$, 
with $\s(a/A)=1$. If $tp_{ACFA}(a/A)\perp (\si(x)=x)$ then $tp(a/A)$ is stable, stably
embedded. In particular, if  $tp_{ACFA}(a/A)$ is stable, stably embedded, then so is $tp(a/A)$. 
\end{prop}

{\it Proof}:\\

 %By \ref{st0}, it suffices to show it when $\s(a/A)=1$.
%Assume that $tp_{ACFA}(a/A)$ is 
%stable stably embedded. 
Suppose that $tp(a/A)$ is not stable, stably embedded; then there is 
$B=acl(B)\supset A$ such that $tp(a/B)$ is not 
stationary, and therefore $tp_{ACFA}(a,Da,D^2a,\ldots /B)$ is not stationary. 

By \ref{st1} 
$tp_{ACFA}(a,Da,D^2a,\ldots /A) \not\perp (\si(x)=x)$. 
Hence, there is some algebraically closed difference field $L$ containing $A$, 
which is linearly disjoint from $acl(Aa)$ over $A$, and an element $b\in Fix(\si)\cap 
(Lacl(Aa))^{alg}, b \not\in L$. Looking at the coefficients of the minimal polynomial of $b$ over $Lacl(Aa)$, 
we may assume that $b\in Lacl(Aa)$.
Let $M=acl(L)$, and chose $(M',L')$ realizing $tp(M,L/A)$ and independent from 
$a$ over $A$. Then $qftp_{ACFA}(L'/Aa)=qftp_{ACFA}(L/Aa)$ and there is $b'\in L'acl(Aa)$ such that $\si(b')=b'$. 
Since $\s(a/L')=1$, we get $a\in acl(L'b')=L(b')_D^{alg}$. This implies that $tp_{ACFA}(a/L')\not\perp (\si(x)=x)$,
 and gives us a contradiction.\\
$\Box$

\begin{rem}\label{st5}
As stated, the result of \ref{st4} is false if one only 
assumes $\s(a/A)<\omega$. The correct formulation in that case is as 
follows:

Assume $\s(a/A)<\omega$ and that $acl_{\sigma}(Aa)$ contains a sequence 
$a_1,\cdots,a_n$ of tuples such that, for all $i\leq n$, working in 
{\it DCFA},  $\s(a_i/Aa_1,\cdots,a_{i-1})=1$.
% or 
%$tp(a_i/Aa_1,\cdots,a_{i-1})$ is almost internal to the set of 
%conjugates of some type of $\s$-rank $1$.
Under these hypotheses, if $tp_{ACFA}(a/A)$ is stable, stably embedded 
then so is $tp(a/A)$.

%Then $tp(a/A)$ is stable 
%stably embedded if and only if $tp_{ACFA}(a/A)$ is stable stably 
%embedded, if and only if for each $i$, 
%$tp_{ACFA}(a_i/Aa_1,\cdots,a_{i-1})\perp (\sigma(x)=x)$. 

\end{rem}
The proof of the following lemma is analogue to the last statement in the proof of \ref{st1}(\ref{st1_2}).t
\begin{lem}\label{dcfdcfa}
Let $a$ be a tuple of a model of {\it DCFA}, and $A$ a subset of that model.
If $tp_{DCF}(a/A)$ is 1-based then $tp(a/A)$ is 1-based.
\end{lem}

%{\it Proof}:\\
%
%Analogue to the last statement in the proof of \ref{st1}.2\\
%Let $A=acl(A) \subset B=acl(B) \subset C=acl(C)$ , and let   $b$ be tuple
 %of realisations of $tp(a/A)$.
% Let us suppose that $acl(Bb) \cap  C=B$.
 %By hypothesis $tp_{DCF}(a/A)$ is
 %$1$-based,  therefore  $b$ is independent from $C$ over $B$ in {\it DCF}.
 %As $A=acl(A)$, we have that for all $n \in \na$, $tp_{DCF}(\si^n(b)/A)$ is
 %$1$-based. By \ref{wgr} $tp_{DCF}(b,\si(b),\ldots, \si^n(b)/B)$ is also
 %$1$-based.
 %Hence, $(b,\si(b),\ldots,\si^n(b))$ is {\it DCF}-independant from $C$ over $B$,
 %for every $n \in  \na$.
 %Then for every finite subset $S$ of 
%$acl(Bb)$, $B(S)$ is linearly disjoint from $C$ over $B$ (that is 
%because 
%every such $S$ is such that $B(S)$ contained in $acl_D(B,b,\si(b), 
%\cdots, 
%\si^n(b))$ for some $n$). Thus by definition of linear disjointness 
%$acl(Bb)$ is linearly disjoint from $C$ over $B$. So
%$b$ is {\it DCFA}-independent from $C$ over $B$.\\
%$\Box$

%The following theorem due to Wagner (\cite{waggr}) gives us conditions for 1-basedness, stability and stable embeddability
%for groups.
Lemmas 2 and 3 of \cite{salinas} and \ref{wgr} imply the following condition for  1-basedness, stability and stable embeddability for groups.

\begin{teo} \label{wgrsec}
Let $1 \longrightarrow G_1   \longrightarrow G_2    \longrightarrow G_3  \longrightarrow  1$ be a
short exact sequence of definable groups in a simple theory. Then $G_2$ is
 stable, stably embedded (resp. 1-based) if and only if $G_1$ and $G_3$ are stable, stably embedded (resp. 1-based).
\end{teo}

\section{Abelian Groups}\label{sec:abelian}

In this section, we study abelian groups defined over some 
subset $K=acl(K)$ of a model $({\mathcal U},\si,D)$ of {\it DCFA}. 
We investigate whether they 
are $1$-based, and whether they are stable, stably embedded. By 4.3 of \cite{rbrank1}, \ref{GIII6} and \ref{wgrsec}  this study may be reduced to the case when the group $H$ is a quantifier-free 
definable subgroup of some commutative algebraic group $G$, and $G$ has no proper (infinite) algebraic subgroup, i.e. $G$ is either $\gr_a$, 
$\gr_m$, or a simple abelian variety $A$. \\\\
{\bf From now on we suppose all the groups are quantifier-free definable}.\\

We study now all three cases for $G$.\\\\
{\bf The additive group}
\begin{prop}\label{adgr}
No infinite  definable subgroup of  $\gr_a^n(\mathcal{U})$ is $1$-based.
\end{prop}

{\it Proof}:\\

Let $H<\gr_a^n$ be a definable infinite group. By 4.4 of \cite{rbrank1},
$H$ is quantifier-free definable and contains a definable subgroup $H_0$ which is definably isomorphic to 
$Fix (\si) \cap {\mathcal C}$. Hence $H$ is not $1$-based.\\
%It is clear that $H$ is defined by linear difference-differential equations, and
%therefore it is a $(Fix \si \cap {\mathcal C})$-vector space.\\
$\Box$\\\\
{\bf The multiplicative group}\\

The logarithmic derivative $lD:\gr_m \to \gr_a$, $x\mapsto Dx/x$ is a group epimorphism with
$Ker(lD)=\gr_m(\mathcal{C})$ (see \cite{manin}). 

Given a polynomial $P(T)= \sum_{i=0}^na_iT^i\in {\mathbb Z}[T]$, we denote by $P(\si)$ the homomorphism defined by
$x \mapsto \sum_{i=0}^na_i\si^i(x)$.

\begin{prop}
Let $H$ be a quantifier-free ${\mathcal L}_{\si,D}$-definable subgroup of $\gr_m$. 
If $lD(H) \neq 0$ then $H$ is not $1$-based. If $lD(H)=0$ then there is a polynomial
$P(T)$ such that $H=Ker(P(\si))$. Then we have that $H$ is $1$-based if and only if $P(T)$ 
is relatively prime to all cyclotomic 
 polynomials $T^m-1$ for all $m\in \na$ 
\end{prop}

{\it Proof}: \\

By \ref{adgr}, if $lD(H) \neq0$ then $H$ is not $1$-based.
If $lD(H)=0$, as $Ker(lD)=\gr_m(\mathcal{C})$, $H$ is ${\mathcal L}_{\si}$-definable in $\mathcal{C}$. 
Hence there is a polynomial 
$P(T)=\sum_{i=0}^na_iT^i \in \mathbb{Z}[T]$ such that $H$ is defined by $\Pi_{i=0}^n\si^i(X^{a_i})=1$.
In {\it ACFA}, $H$ is  1-based, stable, stably embedded if and only if $P(T)$ is relatively prime to all cyclotomic 
 polynomials $T^m-1$ 
for $m \geq 1$ (see \cite{HMM}). By \ref{st1} the same holds for {\it DCFA}.\\
$\Box$\\\\
{\bf Abelian varieties}\\

%This relation has lead to model theoretic proofs of Manin-Mumfurd and
%Mordell Lang conjectures both by E. Hrushovski (\cite{HMM} and \cite{HML}). Recently Pillay found
%more direct proofs to this conjectures, for the case of characteristic zero (\cite{PMM} and \cite{PML}).
%First we mention some facts about abelian varieties in difference and differential fields.
%For a detailed exposition on abelian varieties the reader may consult \cite{langabelian}.                          
\begin{defi}\label{ab1}
An abelian variety is a connected algebraic group $A$ which is complete, that is, for any variety $V$ the
projection $\pi:A \times V \to V$ is a closed map.
\end{defi}

As a consequence of the definition we have that an abelian variety is commutative.

Let $B$ be an algebraic subgroup of an abelian variety $A$. Then $A/B$ is an abelian variety. If in addition $B$
is connected $B$ is an abelian variety.
An abelian variety is called simple if it has no infinite proper abelian subvarieties.
Let $A$ and $B$ be two abelian varieties. Let $f:A \to B$ be a homomorphism. We say that $f$ is an isogeny if
$f$ is surjective and $Ker(f)$ is finite. We say that $A$ and $B$ are isogenous if  there are isogenies
$f:A \to B$ and $g:B \to A$.

\begin{prop}\label{ab5}\textup{({\it ACF}, \cite{langabelian})}
There is no nontrivial algebraic homomorphism from a vector group into an abelian variety.
\end{prop}
Now we mention some properties concerning 1-basedness of abelian varieties in difference and differential fields.

Consider a saturated model $({\mathcal U},\si)$ of {\it ACFA}. In \cite{HMM}, Hrushovski gives a full description of definable subgroups 
of $A({\mathcal U})$ when $A$ is a simple abelian variety defined over ${\mathcal U}$.
 When $A$ is defined over $Fix(\si)$, this description is particularly simple, at least up to commensurability. 
Let $R=End(A)$ (the ring of algebraic endomorphisms of $A$). If $P(T)=\sum_{i=0}^n e_iT^i\in R[T]$,
 define $Ker(P(\si))=\{a\in A({\mathcal U})\mid \sum_{i=0}^n e_i(\si^i(a))=0\}$.

\begin{prop} \textup{({\it ACFA}, \cite{HMM})}\label{ab2} 
Let $A$ be a simple abelian variety defined over $\mathcal U$, and let $B$ be a definable subgroup of 
$A({\mathcal U})$ of finite $\s$-rank. 
\begin{enumerate}
\item If $A$ is not isomorphic to an abelian variety defined over $(Fix(\si))^{alg}$, 
then $B$ is $1$-based and stable, stably embedded. 
\item Assume that $A$ is defined over $Fix(\si)$. 
Then there is $P(T)\in R[T]$ such that $B\cap Ker(P(\si))$ has finite index in $B$ and in $Ker(P(\si))$. 
Then $B$ is $1$-based if and only if the polynomial $P(T)$ is relatively prime to all cyclotomic polynomials 
$T^m-1$, $m\in\na$. If $B$ is $1$-based, then it is also stable, stably embedded.
\end{enumerate} 
\end{prop}
%\begin{prop}\label{ab2}(ACFA, see \cite{HMM})\
%If $A$ is defined over $Fix \si$, then
%$Ker(P(\si))(A)$ is stable,stably embedded and 1-based if and only if $P(T)$ is relatively prime to all cyclotomic
%polynomials $T^n-1$, $n \in \na$.
%\end{prop}
%We make an observation. When we work in characteristic 0, in {\it ACFA} if $tp(a/E)$ is hereditary
%orthogonal to $Fix\si$ then it is stable and stably embedded. This comes from the fact that,
%given $a,F$ and $E' \supset E$ such that $a \dnfo_{E'} F $, $acl_{\si}(E'a)F$ has no finite $\si$-stable extension,
% and this
%implies that $tp_{ACFA}(a/E')$ is stationary.
%From this we have that if a type is 1-based in $\mathcal{C}$ then it is 1-based in $\mathcal{U}$
%(see Chapter 2, section 5).
%\begin{defi}\label{ab3}
%Let $G$ be an algebraic group. By an extension of $G$ by a vector group we mean %an algebraic group $H$
%together with an epimorphism $p: H \to G$ such that $Ker(p)$ is a vector group.
%\end{defi}
%\begin{prop}\label{ab4}
%Let $A$ be an Abelian variety. Then there is a universal extension of $A$ by a vector group. That is,
%there is an extension of $A$ by a vector group $({\hat A},p)$, such that if $(B,q)$ is an extension
%of $A$ by a vector group, then there is $j:{\hat A} \to B$ such that $q = p \circ j$.
%\end{prop}
%As a consequence of \ref{ab4} we have that there is an algebraic section $s:{\hat A} \to \tau({\hat A})$.
We work now in a saturated model $({\mathcal U},D)$ of {\it DCF}. The following is proved in \cite{manin}.

\begin{prop}\label{ab6}
Let $A$ be an abelian variety. Then there is a ${\mathcal L}_D$-definable (canonical) homomorphism
$\mu: A \to \gr_a^n$, for $n=dim(A)$, such that $Ker(\mu)$ has finite Morley rank (a generalization of the notion of algebraic dimension).
\end{prop}
$Ker(\mu)$, is known as the Manin kernel of $A$, we denote it by $A^{\sharp}$. 

\begin{prop}\textup{(Properties of the Manin Kernel, see \cite{manin} for the proofs)}\label{propmanin}\\
Let $A$ and $B$ be abelian varieties. Then
\begin{enumerate}
\item $A^{\sharp}$ is the Kolchin closure of the torsion subgroup $Tor(A)$ of $A$.
\item $(A\times B)^{\sharp}=A^{\sharp}\times B^{\sharp}$, and if $B<A$ then $B\cap A^\#=B^\#$.
%hence if $B<A$, then $B\cap A^{\sharp}=B^{\sharp}$ 
%(since $A$ is isogenous to the direct product of $B$ with an Abelian variety). 
\item A differential isogeny  between
$A^{\sharp}$ and $B^{\sharp}$ is the restriction of an algebraic isogeny from $A$ to $B$.
\end{enumerate}
\end{prop}
We say that an abelian variety descends to the constants if it is isomorphic to an abelian variety defined over the constants.
%\end{defi}

\begin{prop}\label{ab8}\textup{({\it DCF}, see \cite{manin})}
Let $A$ be a simple abelian variety. If $A$ is defined over ${\mathcal C}$, then $A^{\sharp}=A({\mathcal C})$.
If $A$ does not descend to the constants, then $A^{\sharp}$ is strongly minimal and 1-based.
\end{prop}

We now return to {\it DCFA} and fix a saturated model $({\mathcal U,\si,D})$ of {\it DCFA} 
and a simple abelian variety $A$ defined over $K=acl(K) \subset \mathcal{U}$.

Let $H$ be an ${\mathcal L}_{\si,D}$-definable connected subgroup of $A$ defined over 
the difference-differential field $K$ and let $\tilde{H}$ be its $(\si,D)$-Zariski closure.
 Since $H$ is $1$-based if and only if  $\tilde{H}$
is $1$-based (see 4.3 and 4.4 of  \cite{rbrank1}),  we can suppose that $H$ is quantifier-free definable and quantifier-free connected.
%Then there is $k \in \na$ and an ${\mathcal L}_D$-definable subgroup $S$ of
%$A \times \cdots \times A^{{\si}^k}$ such that
%$H=\{x \in A:(x,\si(x), \cdots, \si^k(x)) \in S  \}$.

%If $G=\gr_a^n$, we know that there are no 1-based subgroups of an additive group.
%If $G=\gr_m^n$, then the $D$-definable subgroups of $G$ are defined by linear
%equations in $x,Dx, \cdots,D^mx$ for some $m$, then it is non 1-based because
%it is a $(Fix \si \cap \mathcal{C})$-vector space.\\
Let $\mu:A \to \gr_a^d$ as in \ref{ab6}. If $H \not\subset Ker\mu$ then by \ref{adgr}
$H$ is not $1$-based.

Assume that $H \subset A^{\sharp}$. 
%Then $S \subset A^{\sharp} \times (A^{\sharp})^{\si} \times \cdots \times (A^{\sharp})^{\si^k}$.
We first show a very useful lemma.

\begin{lem}\label{ab20}
Let $H$ be a quantifier-free definable subgroup of $A^{\sharp}$ 
which is quantifier-free connected. Then $H=H'\cap A^{\sharp}$ for some 
quantifier-free ${\mathcal L}_\si$-definable subgroup $H'$ of $A$.
%Let $H$ be a quantifier-free ${\mathcal L}_{\si,D}$-definable subgroup of 
%$A^{\sharp}$. Then there is $k$ and a differential 
%subgroup $S$ of $A^{\sharp}\times\cdots\times (A^{\sharp})^{\si^k}$ such that 
%$H=\{a\in A :(a,\si(a),\ldots,\si^k(a))\in S\}$. We 
%know that $S=U^{\sharp}$ for some algebraic subgroup $U$ of $A\times 
%\cdots\times A^{\si^k}$. By \ref{propmanin}, the results of \cite{HMM} give a full 
%description of the quantifier-free ${\mathcal L}_{\si,D}$-definable 
%subgroups of $A^{\sharp}$. 
\end{lem}

{\it Proof}:\\

%Let $H$ be a quantifier-free definable  subgroup of $A^{\sharp}$ which is connected for the $(\si,D)$-topology. 
%Then there is an integer
%$k$ and  a connected algebraic subgroup $S$ of $A\times A^\si\times 
%A^{\si^k}$ such that $H=\{a\in A^{\sharp} : (a,\si(a),\cdots,\si^k(a)\in S\}$.
Our hypotheses imply that there is an integer $k$ and a differential 
 subgroup $S$ of $A\times A^\si\times \cdots \times A^{\si^k}$ such that 
$H=\{a\in A: (a,\si(a),\cdots,\si^k(a))\in S\}$. By \ref{propmanin}.2, replacing $S$ 
by its Zariski closure $\bar S$ we get $H=\{a\in A^{\sharp}: 
(a,\si(a),\cdots,\si^k(a))\in \bar S\}$. 
Thus $H=H'\cap A^{\sharp}$, with $H'=\{a\in A : (a,\si(a),\cdots,\si^k(a)\in
\bar{S}\}$.\\
$\Box$

Let us state an immediate consequence of \ref{ab20} : 

\begin{cor}\label{ab9}
If for all $k \in \na$, $A$ and $A^{\si^k}$ are not isogenous,
then $\s(A^{\sharp})=1$.
\end{cor}
%{\it Proof}:\\
%Otherwise, $A^{\sharp}$ contains an infinite definable proper subgroup; as in \cite{HMM}, we get a $D$-definable 
%isogeny $f:A^{\sharp} \to (A^{\sharp})^{\si^k}$. By \ref{ab6.1} this gives us an algebraic isogeny $f':A \to A ^{\si^k}$.
%
%$\Box$
%\begin{rem}
%With the help of \ref{ab6.1} and \ref{J313}, we can prove that if $A$ is not isomorphic to an Abelian variety defined over 
%$Fix \si^k$ then every finite-dimensional subgroup of $A^{\sharp}$ is stable (the proof is the same as in \cite{HMM}).
%\end{rem}      
%\begin{prop}\label{ab7}
%Let $A$ be a simple Abelian variety. Then $A^{\sharp}$ has no infinite $D$-definable subgroups.
%\end{prop}
\
{\bf Case 1}: $A$ is isomorphic to a simple abelian variety $A'$ defined over $\mathcal{C}$.

We can suppose that $A$ is defined over $\mathcal{C}$. Then, by \ref{ab8},
$A^{\sharp}=A(\mathcal{C})$.
% Hence $S \subset A(\mathcal{C}) \times A^{\si}(\mathcal{C}) \times \cdots \times A^{\si^k}(\mathcal{C})$. 
Hence, by \ref{st1},  $H$ is $1$-based for {\it DCFA} if and only if it is $1$-based for $ACFA$;
and in that case, by \ref{st3}, it will also be stable, stably embedded 

 If $H=A(\mathcal{C})$ then we know that $H$ is not 1-based in {\it ACFA}.
%and
% this implies that $H$ is 1-based for {\it DCFA}: From the way independence is defined we have that if $tp_{DCF}(a)$ then $tp(a)$ is 1-based, similarly, if
 %$tp(a)$ is not 1-based then one of the types $tp_{DCF}(D^ka/a \cdots D^{k-1}a)$
 %is not 1-based.

If $H$ is a proper subgroup of $A(\mathcal{C})$, 
\ref{ab2} gives a precise 
description of 
that case.
%, hence, by \ref{ab2}, $H$ is not 1-based if and
%only if some type in the semi-minimal analysis of $S$ is non orthogonal to $Fix \si \cap {\mathcal C}$ (see \cite{HMM}),
%and again the results of that paper give a precise description of this case.
\\

{\bf Case 2}: $A$ does not descend to $\mathcal{C}$.

Then, by \cite{manin}, section 5, $A^{\sharp}$ is strongly minimal and 1-based for {\it DCF}. By \ref{dcfdcfa} it is 1-based
for {\it DCFA}.

%
%
%If $H \neq A^{\sharp}$ then $dim(H) < \infty$. Then if $A$ is not isomorphic to an Abelian veriety defined over
%$Fix \si^k$ then $H$ is stable and stably embedded.
%Then $H$ is not 1-based if and
%only if some type in the semi-minimal analysis is nonorthogonal to $Fix\si \cap \mathcal{C}$.
%If $H=A^{\sharp}$, as $A^{\sharp}$ is 1-based for {\it DCF} then it is 1-based for {\it DCFA}.

%\begin{teo}
%$\s(A^{\sharp})=1$.
%\end{teo}
%
%{\it Proof}:\\
%
%Suppose that there is an isogeny between $A$ and $A^{\si^k}$ for some $k$, then as in \cite{HMM},
%there is an $D$-definable isogeny
%$f:A^{\sharp} \to (A^{\sharp})^{\si^k}$. Without lose of generality we may
%suppose that $k=1$.
%
%So $f$ gives us a $D$-definable subgroup $H$ of
%$A^{\sharp} \times (A^{\sharp})^{\si}$. Then we have that there is $m \in \na$ and an algebraic
%subgroup $H'$ of $\tau_m(A) \times \tau_m(A^{\si})$ projecting onto $A$ and onto $A^{\si}$, such that
%$H=\{x \in A: (x,Dx,\cdots,D^m x) \in H'\}$. 

We will now investigate when $H$ is stable, stably embedded. By $1$-basedness and quantifier-free 
$\omega$-stability, we know that if $X\subset A^{\sharp}$ is quantifier-free 
definable, then $X$ is a Boolean combination of cosets of 
quantifier-free definable subgroups of $A^{\sharp}$.  

Assume first that $H\neq A^{\sharp}$, and let $a$ be a generic of $H$ over $K$. 
Then $H$ is finite-dimensional, and therefore $\s(H)<\omega$.
As $H$ is $1$-based, there is an increasing sequence 
of subgroups $H_i$ of $H$ with  $\s(H_{i+1}/H_i)=1$.

By \ref{ab20}, we may assume that $H_i=U_i\cap A^{\sharp}$ for some quantifier-free 
${\mathcal L}_\si$-definable subgroups $U_i$ of $A$. Note that \ref{ab20} also 
implies that each quotient $U_{i+1}/U_i$ is $c$-minimal (i.e., all
quantifier-free definable ${\mathcal L}_\si$-definable subgroups are 
either 
finite or of finite index). Furthermore, by elimination of imaginaries 
in 
{\it ACFA}, $acl_{\si}(Ka)$ contains tuples $a_i$ coding the cosets $a+U_i$. 
Hence $tp(a/K)$ satisfies the conditions of \ref{st5} and we obtain that 
if $tp_{ACFA}(a/K)$ is stable, stably embedded then so is $tp(a/K)$. 

For the other direction, observe that if $tp_{ACFA}(a/K)$ is not 
stable, stably embedded, then for some $i$, the generic {\it ACFA}-type of 
$U_{i+1}/U_i$ is non-orthogonal to $\si(x)=x$, and there is a (${\mathcal 
L}_\si$)-definable morphism $\psi$  with finite kernel $U_{i+1}/U_i \to 
 B(Fix (\si^k))$ 
for some $k$ and abelian variety $B$ (see \cite{HMM}). 
But, returning to {\it DCFA}, no non-algebraic type realized in $Fix( \si^k)$ can 
be stable, stably embedded, since for instance the formula 
$\varphi(x,y)=\exists z\ z^2=x+y \ \land \ \si(z)=z$ is not definable (\ref{st1},3). 
This proves the other implication.

\smallskip

%Assume first that $H\neq A^{\sharp}$. Then $H$ is finite-dimensional, and therefore has finite $\s$-rank. 
%By the discussion above,  we may assume that $\s(H)=1$. Let $a,m$ be as above, and consider 
%$tp_{ACFA}(a,Da,\cdots,D^ma/A)$. By \ref{st2}, $tp(a/K)$ will be stable and 
%stably embedded if and only if it is stationary. Assume that 
%$tp_{ACFA}(a,Da,\cdots,D^ma)$ is not stationary. By 4.11 of \cite{salinas}, this 
%implies that there is an algebraically closed  difference field $L$ 
%containing $K$ and linearly disjoint from $K(a,Da,\cdots,D^ma)_\si$ 
%over $K$, and an element $b\in Fix\si$ such that $b\in 
%acl_{\si}(L,a,\cdots,D^ma)$. Looking at the coeeficients of the minimal 
%polynomial of $b$ over $L(a,\cdots,D^ma)_{\si}$, we may assume that 
%$b\in L(a,\cdots,D^ma)_{\si}$. Hence, we may also assume that 
%$L=acl(L)$. Since $\s(a/L)=1$, this implies that $a\in 
%acl(L,b)=L(b)_D^{alg}$. Since all $D^ib\in Fix\si$, we get 
%$tp_{ACFA}(a/L)\not\perp \si(x)=x$. 
Thus we have shown: 

\smallskip 
If $H$ is finite dimensional, then $tp(a/K)$ is stable, stably embedded 
if and only if $tp_{ACFA}(a/K)$ is stable, stably embedded. 

%The results in \cite{HMM} then give a complete description of that case. 
Using \ref{ab20}, \ref{ab2} gives us a full description of that 
case.

In particular, we then have that if $H$ is not stable, 
stably embedded, then $A$ is isomorphic to an abelian 
variety defined over $Fix(\si^k)$ for some $k$. 

\smallskip 
Let us now assume that $H=A^{\sharp}$. Let $a$ be a generic of $H$ over $K$. Then 
$tp_{ACFA}(a,\cdots,D^ma/K)$ is the generic type of an algebraic 
variety $V$, and is therefore stationary (by 2.11 of \cite{salinas}). Thus, using the finite dimensional case, if 
$A$ is not isomorphic to an abelian variety defined over 
$(Fix(\si))^{alg}$, then $H$ is stable, stably embedded. If $A$ is 
isomorphic to a variety $B$ defined over $Fix( \si^k)$, via an 
isomorphism $\psi$, then the subgroup $\psi^{-1}(Ker (\si^k -1))\cap A^{\sharp}$ is 
not stable, stably embedded. 
%Indeed, we may assume that  $k$ is large 
%enough so that $B[2]\subset B(Fix(\si^k))$. Let $b\in B$, $\si^k(b)=b$, $b$ 
%differentially transcendental over $K$. Then 
%$tp(\psi(a)/K)$ has several distinct non-forking extensions to $acl(Kb)$: 
%those containing 
%the formula $\exists y \in B\, 2y=x+b \land \si^k(y)=y$, and those 
%containing its negation. 

We summarize the results obtained:

\begin{teo} \label{absum} 
Let $A$ be a simple abelian variety, and let $H$ be a 
quantifier-free definable subgroup of $A({\mathcal U})$ defined over $K=acl(K)$. If 
$H\not\subset A^{\sharp}({\mathcal U})$, then $H$ is not 1-based. Assume 
now that $H\subset A^{\sharp}({\mathcal U})$, and let $a$ be a generic of 
$H$ over $K$. Then 
\begin{enumerate}

\item If $A$ is defined over the field  $\mathcal C$ of constants, then 
$H$ is $1$-based if and only if it is stable, stably embedded, if and 
only if $tp_{ACFA}(a/K)$ is hereditarily orthogonal to 
$(\si(x)=x)$. The results in \cite{HMM} yield a complete description of 
the subgroups $H$ which are not $1$-based. 

\item If $A$ does not descend to the field  $\mathcal C$ of constants, 
then $H$ is $1$-based. Moreover 

\begin{enumerate}
\item If $A$ is not isomorphic to an abelian variety defined over 
$Fix( \si^k)$ for some $k$, then $H$ is stable, stably embedded. 

\item Assume that $A$ is defined over $Fix(\si)$. Then $H$ is stable, 
stably embedded if and only $tp_{ACFA}(a/K)$ is stable, stably 
embedded. Again, the results in  \cite{HMM} give a full description of 
this case. 
\end{enumerate}
\end{enumerate}
\end{teo}

\def\rasp{\leavevmode\raise.45ex\hbox{$\rhook$}} \def\cprime{$'$}

\bibliographystyle{plain}
\end{document}